\documentclass[smallextended,numbook,runningheads]{svjour3}     
\smartqed  
\usepackage{graphicx}
\DeclareGraphicsExtensions{.ps,.eps,.png}
\usepackage{epsfig}
\usepackage{mathptmx}      
\usepackage{verbatim}
\usepackage{amsmath}

\usepackage{subfigure}

\newcommand\asp{\mathop{\rm asp}}
\newcommand\minalt{\mathop{\rm minalt}}

\journalname{BIT}
\newcommand\bz{{\bf 0}}
\newcommand\commentout[1]{\relax}
\newcommand\R{{\bf R}}

\newcommand\x{{\bf x}}
\newcommand\eref[1]{$(\ref{#1})$}

\newcommand\tabref[1]{Table \ref{#1}}

\begin{document}

\title{Analysis of and workarounds for element reversal for a finite 
element-based algorithm for warping triangular and tetrahedral meshes
\thanks{The majority of this work was performed while the first author
was a member of the Center for Applied Mathematics at Cornell University.
The work of the first author was supported by the National Physical Science
Consortium, Sandia National Laboratories, and Cornell University.
The work of both authors was supported in part by NSF grant ACI-0085969.}\footnote[1]{Paper accepted for 
publication in {\emph{BIT Numerical Mathematics}}.  The final publication is available at 
www.springerlink.com (DOI: 10.1007/s10543-010-0283-3).}
}

\titlerunning{Analysis of and Workarounds for Element Reversal in Tetrahedral Mesh Warping}

\author{Suzanne~M.~Shontz         \and
        Stephen~A.~Vavasis 
}

\authorrunning{Shontz and Vavasis} 

\institute{Suzanne~M.~Shontz \at
              Department of Computer Science and Engineering, The Pennsylvania State University, 
              University \\
              Park, PA  16802  United States of America\\
              Tel.: +1-814-865-0193\\
              Fax: +1-814-865-3176\\
              \email{shontz@cse.psu.edu}           
           \and
           Stephen~A.~Vavasis \at
              Department of Combinatorics and Optimization, University of Waterloo, Waterloo, \\ 
              Ontario, N2L 3G1, Canada\\
              \email{vavasis@math.uwaterloo.ca}}

\date{Received: date / Accepted: date}

\maketitle

\begin{abstract}
We consider an algorithm called FEMWARP for warping triangular 
and tetrahedral finite element meshes that computes the warping 
using the finite element method itself.  The algorithm takes as input
a two- or three-dimensional domain defined by a  boundary mesh (segments in one dimension or
triangles in two dimensions) that has a volume mesh (triangles in two dimensions or tetrahedra in three
dimensions) in its interior.
It also takes as input
a prescribed movement of the
boundary mesh.  It computes as output
updated positions of
the vertices of the volume mesh.  The first step of the algorithm is to determine
from the initial mesh
a set of local weights for each interior vertex that describes each
interior vertex in terms of the positions of its neighbors.  These weights
are computed
using a finite element stiffness matrix.
After a
boundary transformation is applied, a linear system of equations based
upon the weights is solved to determine the final positions of the
interior vertices.

The FEMWARP algorithm has been considered in the previous literature (e.g., 
in a 2001 paper by Baker).
FEMWARP has been succesful in computing deformed meshes for certain
applications.
However, sometimes FEMWARP reverses
elements; this is our main concern in this paper.
We analyze the causes for this undesirable behavior and propose 
several techniques to make the method more robust against 
reversals.  The most successful of the proposed methods includes combining FEMWARP 
with an optimization-based untangler.
\end{abstract}

\section{Introduction}
\label{sec:intro}
There are numerous applications in science and engineering (and other domains as well) for which the domain of 
interest deforms as a function of time.  These applications include structural 
anatomical remodeling associated with heart failure~\cite{helm}, traumatic brain injuries~\cite{darvish}, 
metal forming~\cite{lee}, and fluid flow applications~\cite{ishikawa} (to name just a few).  For such 
applications, the mesh must be updated at each time step in response to the deforming domain boundary, 
thus resulting in potentially drastically varying mesh quality from step to step.  It is well-known that poor 
quality elements affect the stability, convergence, and accuracy of finite element and other solvers 
because they result in poorly-conditioned stiffness matrices and poor solution 
approximation~\cite{shew_quality}.  Hence, it is important for the updated meshes to be of reasonable
quality, as well. 

The specific mesh-updating problem on deforming domains which we study is as 
follows.  Here we follow the description in~\cite{knupp_dd}.  Suppose $\Omega$ is an initial undeformed domain 
and $M$ is an original mesh on $\Omega$.  Further, suppose that $M$ is of good quality.  Then, suppose the 
domain deforms due to a change in shape, orientation, volume, etc.  Let $\hat{\Omega}$ be the 
deformed domain.  The goal is to create a mesh $\hat{M}$ on $\hat{\Omega}$ that is of reasonable quality.  
In addition, it is often desired that $\hat{M}$ be `similar' to $M$ in some sense; e.g., we desire 
that the two meshes have the same topology. Similar meshes are often desired between successive 
time steps of a physical simulation or between successive iterations of an optimal-design procedure 
so that the solution varies smoothly between successive deformations.

If one were given an onto map $X(\Omega,t)$ from $\Omega$ to $\hat{\Omega}$ at a specified time $t$, 
then $\hat{M}$ could be determined by evaluating $X(\Omega,t)$ at the vertices of $M.$ However, for 
many applications, such a map is not known.  In this situation, $\hat{M}$ could be created via an 
automatic mesh generator. However, this is a computationally intensive process, and the resulting 
mesh will likely not be similar to $M$.  Thus, $\hat{M}$ should be created using a mesh-update 
procedure.  Even if the boundary map between $\Omega$ and $\hat{\Omega}$ is given, the mesh-updating 
problem is not easy due to the similarity requirement.

Several partial differential equation-based (i.e., PDE-based) approaches for solving the mesh-updating problem 
in response to a deforming domain boundary, under the assumption that the boundary map, $X(\Omega,t)$ from 
$\Omega$ onto $\hat{\Omega}$ at time $t$, is specified, have been designed.  Work has focused on the 
development of spring model approaches based on Laplace's equation, variable diffusion, and biharmonic PDEs for vertex 
movement~\cite{helenbrook2003,branets2005,shontz2003,stein2003}. Other research has focused on the 
development of elasticity-based approaches~\cite{stein2004,tezduyar1992,ShontzThesis}.  Many 
existing serial and parallel mesh-updating 
methods~\cite{baker,li2001,cardoze1,cardoze2,selwood96parallel,selwood97parallel,shephard03automated,seol06efficient,antaki00parallel} 
combine vertex movement with other techniques which alter the mesh topology and violate the 
similarity requirement.

Another important limitation of most mesh-updating approaches is that they can reverse elements during the mesh 
updating procedure during a given timestep.  ``Reversal'' means that the element changes orientation.  In two 
dimensions, this means that the vertices of an element are clockwise when they ought to be counterclockwise, 
and in three dimensions it means that they violate the right-hand rule.  Meshes with reversed elements yield 
physically invalid solutions (e.g., when such meshes are used in conjunction with the finite element method in 
order to solve a partial PDE).  Hence, it is crucial that the mesh warping procedure does not reverse any 
elements when deforming the mesh.

An optimization-based approach to the mesh-updating problem, based on the target matrix 
paradigm~\cite{knupp_tmp}, recently appeared in~\cite{knupp_dd}. The approach uses mesh 
optimization to create a mesh similar to and having the same topology as $M$.  However, the 
method is computationally expensive and does not guarantee the prevention of element reversal in 
the deformed mesh.  


Another mesh-updating approach, called FEMWARP, computes the warping itself based on the   
finite element method.  FEMWARP is equivalent to a weighted version of Laplacian smoothing,
i.e., a homogeneous Poisson equation is solved for each interior vertex coordinate with Dirichlet
boundary conditions given by the new vertex coordinates of the boundary.
FEMWARP has been considered in the previous literature, e.g., by Baker~\cite{baker},
although he rejected it in favor of a method based on linear elasticity.  It is shown, however, in
\cite{ShontzThesis} that there appear to be few advantages of linear elasticity over FEMWARP;
whereas, a significant disadvantage is that the linear elasticity matrix problem is three times
larger (in 3D) than FEMWARP's (although FEMWARP must solve the smaller linear system three times).

FEMWARP is described in Section~\ref{sec:femwarp}. There are two main advantages to the FEMWARP 
algorithm. First, if a continuous deformation of the boundary is given, then FEMWARP is valid for 
computing the resulting trajectory that specifies the movement of the interior vertices.  In addition, 
these trajectories will be continuous due to the similarity requirement.  This is vital for some 
applications (e.g., in biomechanics~\cite{sigal}, automotive design~\cite{donders}, computer 
animation~\cite{hu}, and clothing design~\cite{luo}), where continuity of motion is required.  A second 
big advantage is that sparse matrix algorithms may be used to solve (\ref{eqn:lastsystem_framework}).  
The sparsity structure is apparent, since, on average, an interior vertex has six neighbors in 2D, 
whereas a typical 2D mesh may have thousands of vertices.

FEMWARP is exact for affine boundary mappings as proven in Section~\ref{sec:affine}.  However, the principal 
failure mode for FEMWARP is element reversal during the mesh updating procedure during a given 
timestep~\cite{ShontzThesis}.  Element reversal is our main concern in this paper.  The causes of element 
reversal are covered in Section~\ref{sec:inversion}.  Techniques to prevent some reversals, including 
small-step FEMWARP, mesh refinement, and the use of another mapping to compute the weights used to 
determine the warping, are also covered in that section.  Another technique to avoid 
reversals based on the Opt-MS mesh untangler is covered in Section~\ref{sec:untangling}. In 
Section~\ref{sec:test3d}, we test our algorithms on several types of deformations of three-dimensional meshes. 
Concluding remarks are presented in Section~\ref{sec:conclusions}.

\section{The FEMWARP algorithm}
\label{sec:femwarp}
In this section, we describe the three-step FEMWARP algorithm (see, e.g., ~\cite{baker}).
The first step of the FEMWARP algorithm is to express the coordinates
of each interior vertex of the initial mesh as a linear combination of
its neighbors. Let a triangular or tetrahedral mesh, $M$, be given for the
domain $\Omega$ in two or three dimensions.
Let $b$ and
$m$ represent the numbers of boundary and interior vertices, respectively.
Form the $(m+b)\times(m+b)$ stiffness
matrix $A$
based on piecewise linear finite elements defined on the initial mesh
for the boundary value problem
$$\bigtriangleup u =  0 \quad\mbox{on $\Omega$}$$
with $u=u_0 \ \mathrm{on \ \partial \Omega}.$  
Because we only keep the relevant matrix, any 
$u_0$ may be chosen.
It is well-known \cite{johnsonbook} that this matrix is determined as
follows.  Let $\phi_i$ be the 
continuous piecewise linear function (where the
pieces of linearity are given by the triangulation) such that $\phi_i(x_i)=1$,
where $x_i$ is the $i^{th}$ vertex of the mesh, and $\phi_i(x_j)=0$, where
$x_j$ is any other vertex in the mesh ($j\ne i$).  Define for
each $i,j=1,\ldots,m+b$
$$A(i,j)=\int_\Omega \nabla \phi_i\cdot\nabla\phi_j.$$
This matrix will be sparse and symmetric positive semidefinite.
Its nonzero entries correspond to pairs of neighboring vertices in the
mesh.

Next, let $A_I$ denote the $m\times m$ submatrix of $A$ whose
rows and columns
are indexed by interior vertices, and let $A_B$ denote the $m\times b$
submatrix of $A$ whose rows are indexed by interior vertices and whose
columns are indexed by boundary vertices.
Let $x$ be the $(m+b)$ vector consisting of $x$-coordinates of 
the vertices of the initial mesh, where we assume that the interior vertices
are numbered first.
Then it follows from well-known theory
that $[A_I,A_B]x=0$  because any linear function of the coordinates
is in the null-space of the discretized Laplacian operator.  
For the same reason, a similar identity holds for the $y$- and 
$z$-coordinates.  An equivalent way to write this equation is
\begin{equation}
A_Ix_I =-A_Bx_B.
\label{eq:aiab}
\end{equation}

If we divide each row of $[A_I,A_B]$ by the diagonal
element in that row, we obtain a linear system whose diagonal
entries are $1$'s and whose 
row sums are $0$'s.  This means that the
$[A_I,A_B]$, thus scaled, expresses each interior vertex coordinate 
as an affine combination of the neighboring vertex coordinates.  
The sign of these weights is important.  In
particular, in 2D, if the boundary of the original mesh is convex
and the weights are nonnegative, it can be shown that there is
no element reversal (see Theorem 4.1,~\cite{floater1}).  These weights
are nonnegative if and only if the two angles opposite to each mesh
edge sum to at most $\pi$~\cite{floater1}.

The formation of $A_I$ and $A_B$ is the first step of the FEMWARP method.
Consider now the  application of
a user-supplied transformation to the boundary 
of the mesh.  We denote the new positions of the boundary vertices by
$[\hat{x}_B, \hat{y}_B]$ in two dimensions
or $[\hat{x}_B, \hat{y}_B, \hat{z}_B]$ in three.

The final step is to solve a linear system of equations 
similar to (\ref{eq:aiab})
for the new 
coordinates of the interior vertices of $\hat{M}$ on $\hat{\Omega}$. 
We solve (\ref{eqn:lastsystem_framework}) for $[\hat{x}_I,
\hat{y}_I]:$
\begin{equation}
A_I [\hat{x}_I, \hat{y}_I] = -A_B [\hat{x}_B, \hat{y}_B].
\label{eqn:lastsystem_framework}
\end{equation}
or the analog in three dimensions.  Due to the similarity
requirement, the mesh topology of $\hat{M}$ is the same
as that of $M$.  Hence $\hat{M}$ is fully specified after
solving (\ref{eqn:lastsystem_framework}).  This concludes the description of FEMWARP.

\section{FEMWARP is exact for affine mappings}
\label{sec:affine}
One useful property of the FEMWARP algorithm
is that the method is {\em exact} for affine transformations.
Let us state this as a lemma.  The lemma is stated for the
two-dimensional case, and it extends in the obvious way to three
dimensions.
\begin{lemma}
Let $A_B$ and $A_I$ be generated using FEMWARP.
Then $A_I$ is nonsingular based upon well-known finite element 
theory~\cite{johnsonbook}.  
Let $[\hat x_B,\hat y_B]$ be
the user-specified deformed coordinates of the boundary.  
Suppose there exists a $2\times 2$ nonsingular matrix $L$ and
$2$-vector $v$ such that for each $j\in B$, 
$$\left(
\begin{array}{c}
\hat x_j \\
\hat y_j
\end{array}
\right)=
L
\left(
\begin{array}{c}
 x_j \\ 
 y_j
\end{array}
\right)+v.$$
Let $[\hat x_I,\hat y_I]$ be the deformed
interior coordinates computed by the method.  Then for each $i\in I$,
$$\left(
\begin{array}{c}
\hat x_i \\
\hat y_i
\end{array}
\right)=
L
\left(
\begin{array}{c}
 x_i \\ 
 y_i
\end{array}
\right)+v.$$
\label{thm:affine}
\end{lemma}

\begin{proof}
The positions of
the interior vertices in the deformed mesh
are given by
\begin{equation}
[\hat x_I, \hat y_I] = -A_I^{-1}A_B([x_B, y_B]L^T+e_B v^{T})
\label{eqn:unique}
\end{equation}
where, as above, $\hat x_I,\hat y_I$ are column vectors composed of the $x$-
and $y$-coordinates of the interior vertices respectively and
$x_B$, $y_B$ are the corresponding vectors for boundary vertices,
and finally $e_{B}$ is vector of all $1$'s of length $b$.

In order to show that affine mappings yield exact results with any
algorithm within the framework, we
want to show that (\ref{eqn:unique}) is the same as:
\begin{equation}
[\hat x_I,\hat y_I] = [ x_I, y_I]L^T+e_I v^{T}.
\label{eqn:guess}
\end{equation}

Observe that the equivalence of (\ref{eqn:unique}) and 
(\ref{eqn:guess}) would follow immediately from:
\begin{equation}
A_I([x_I,y_I]L^T+e_I v^{T}) = -A_B([x_B,y_B]L^T+e_B v^{T}).
\label{eqn:truth}
\end{equation}

Thus, it remains to check that (\ref{eqn:truth}) holds.

Because the weights for each interior vertex sum to $1$,
$A_I e_I + A_B e_B = \bz$ as noted above.
Hence $(A_I e_I + A_B e_B)v^{T} = \bz.$  
Also, because $[x_I, y_I]$ and $[x_B, y_B]$
denote the original positions of the vertices, we know that $A_I[x_I, 
y_I]+A_B[x_B, y_B] = \bz$.  So
$(A_I [x_I, y_I]+A_B[x_B, y_B])L^T = \bz$.  

Putting these together, we see that
\begin{equation}
(A_I [x_I,y_I]+A_B [x_B,y_B])L^T + (A_I e_I + A_B e_I) v^{T} = \bz.
\label{eqn:good}
\end{equation}

Therefore, (\ref{eqn:truth}) holds, and the lemma is proven.
\end{proof}

\section{Element reversal and small-step FEMWARP}
\label{sec:inversion}
Sometimes FEMWARP is successful in yielding valid deformed meshes, i.e.,
those without any reversed elements.  For example, in~\cite{ShontzThesis}, the 
first author was successful in using FEMWARP to generate a sequence of
deforming meshes for use in studying the beating canine heart
from an atrial pacing experiment.  However, sometimes FEMWARP 
fails to yield a valid triangulation because it reverses elements in the
resulting deformed mesh.  

\begin{figure}[ht]
\begin{center}
\caption{Translation and rotation of annulus mesh:  (a) original mesh;
(b)-(d) meshes obtained by translating the inner circle outwards and
rotating it counterclockwise by $(1,0^{o}), (2,20^{o}), (3,40^{o})$,
and $(3.5,50^{o})$, respectively.  Meshes (a) and (b) are valid meshes;
mesh (c) is also valid but contains elements that are near reversal;
mesh (d) contains reversed elements and is invalid.}
\label{fig:2d_mappings}
\vbox{
\hbox{
\hspace{.5in}
  \subfigure[$r=1, \theta = 0^{o}$]{
    \includegraphics[height=4cm,width=4cm]{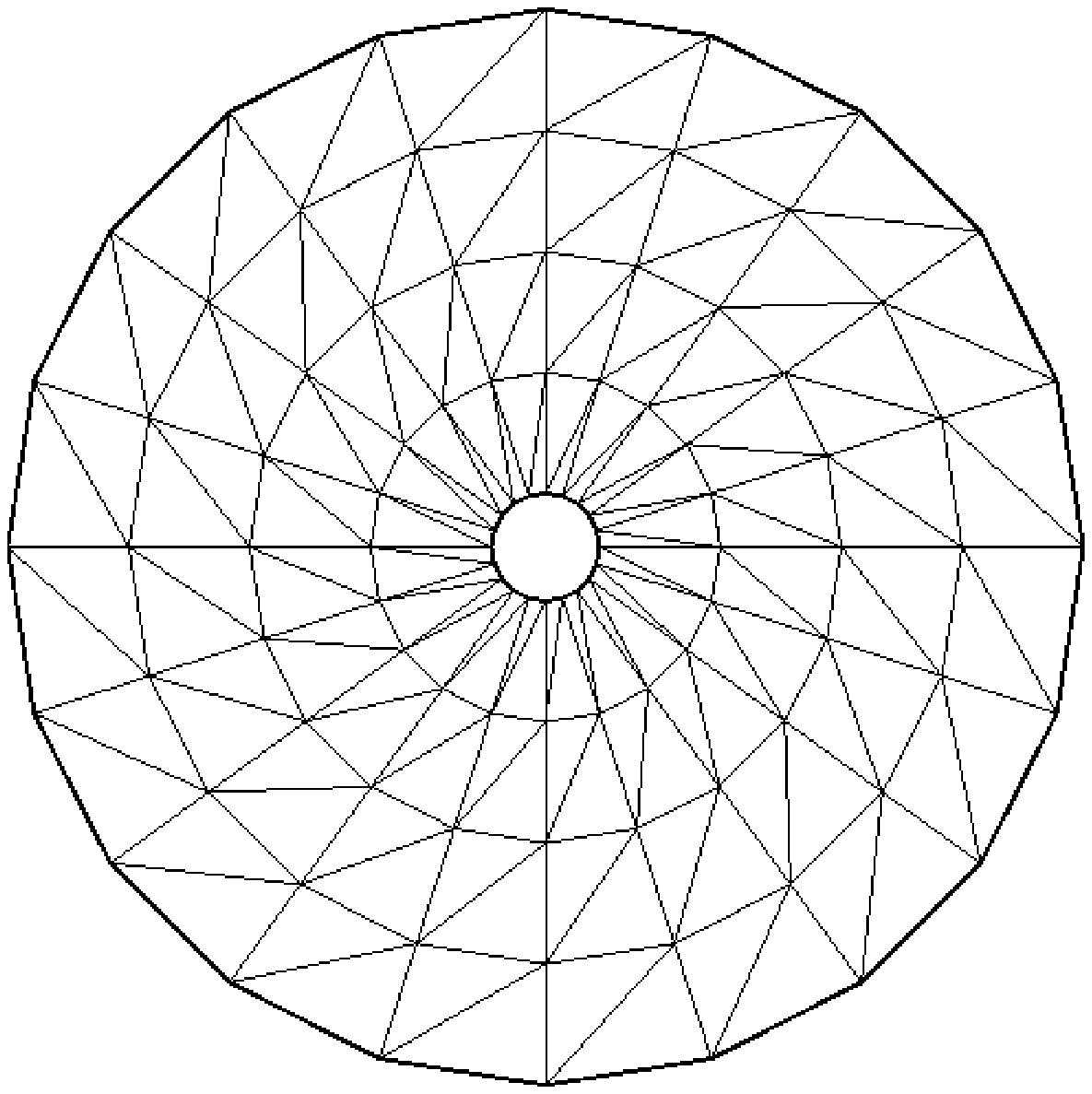}
    \label{fig:annulustheta0}
  }
  \subfigure[$r=2, \theta = 20^{o}$]{
    \includegraphics[height=4cm,width=4cm]{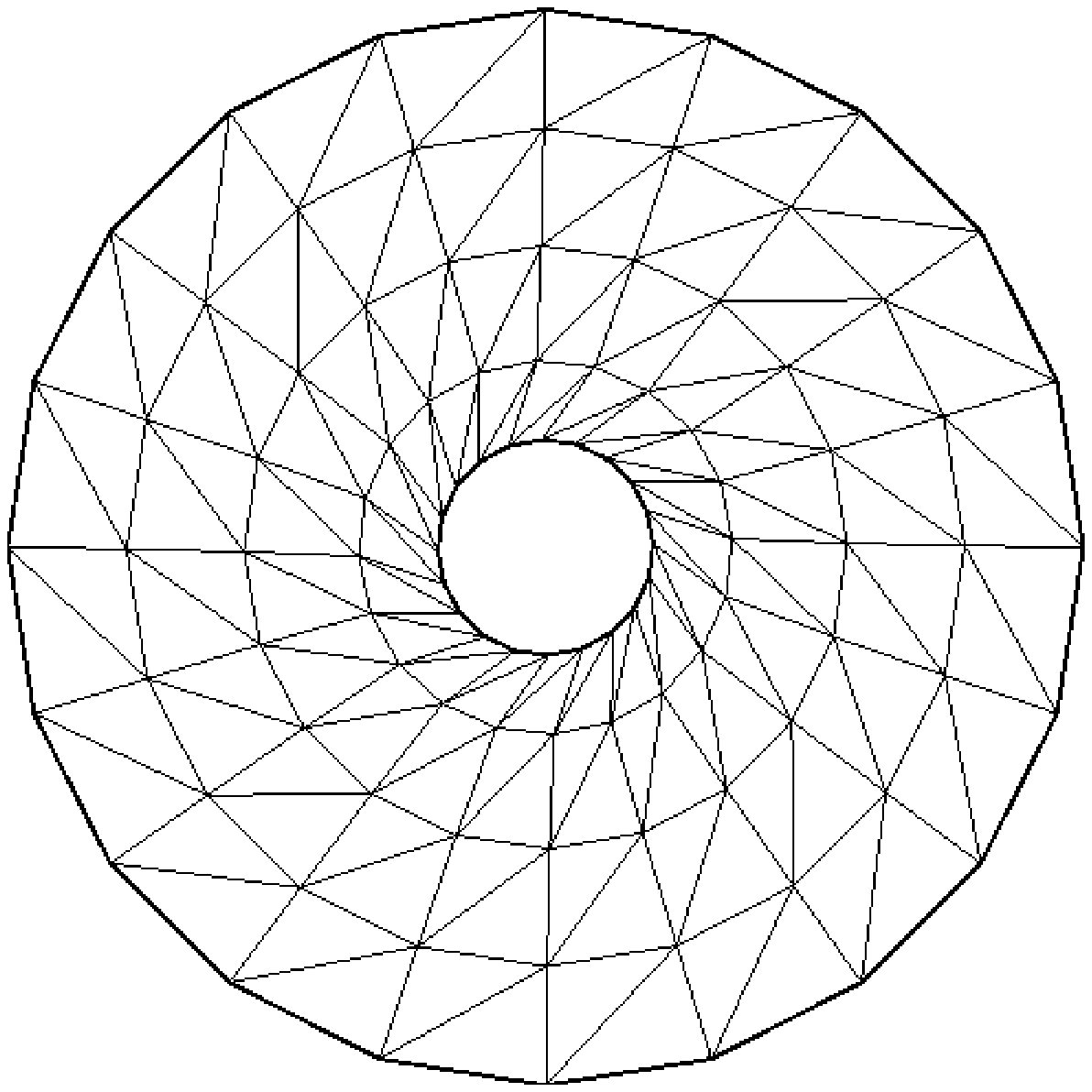}
    \label{fig:annulus_theta20}
  }
}
\hbox{
\hspace{.5in}
\subfigure[$r=3, \theta = 40^{o}$]{
\includegraphics[height=4cm,width=4cm]{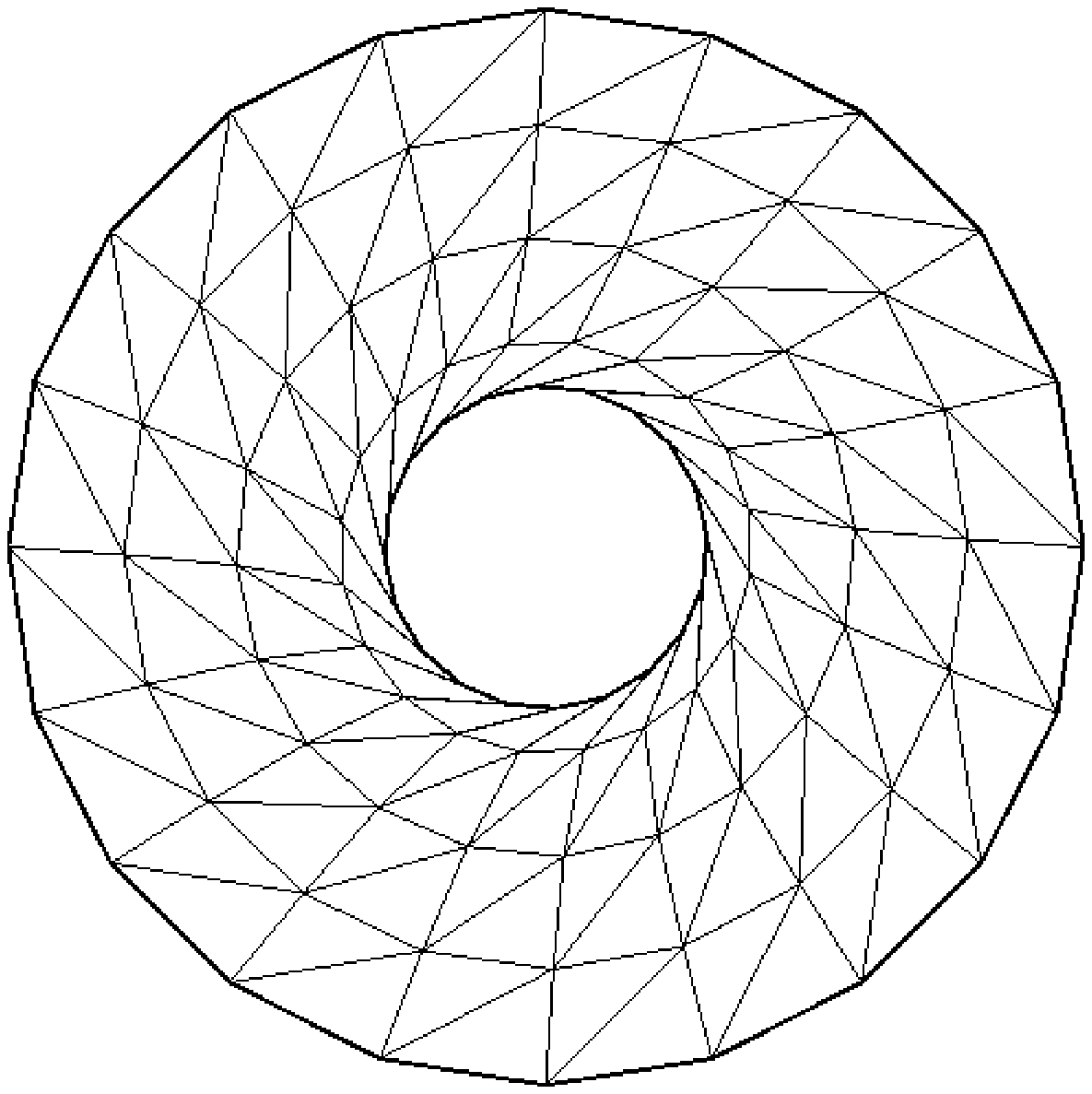}
\label{fig:annulustheta40}
}
\subfigure[$r=3.5, \theta = 50^{o}$]{
\includegraphics[height=4cm,width=4cm]{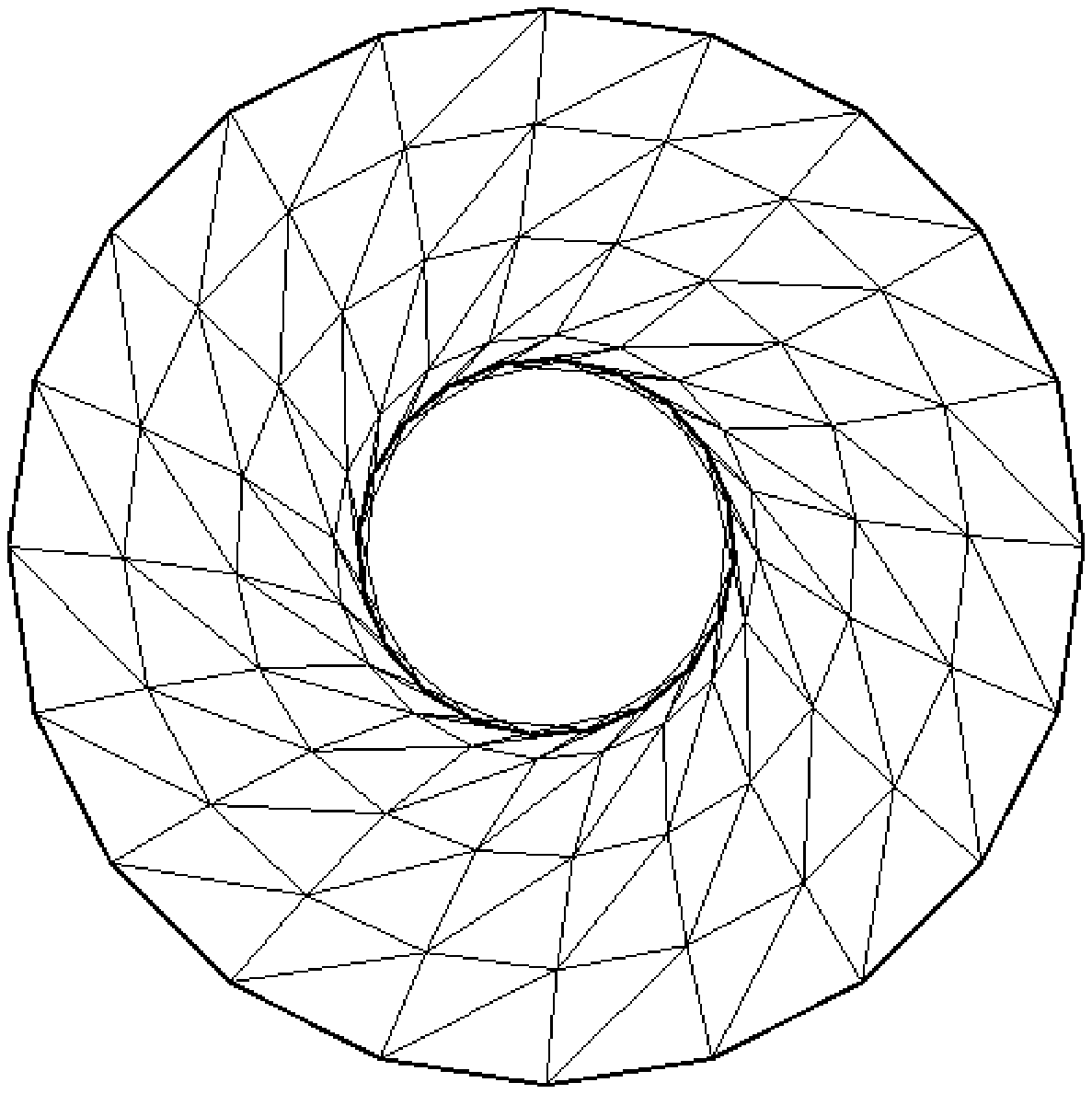}
\label{fig:annulustheta50}
}
}
}
\end{center}
\end{figure}

An example of the valid and invalid deformed meshes which can result using 
FEMWARP is shown in Figure~\ref{fig:2d_mappings}.  For this example,
an annulus mesh (shown in Figure~\ref{fig:2d_mappings}(a)) is deformed.
The annulus mesh in this example is composed of four, equally-spaced
concentric rings of triangles for a total of 160 triangles.  Its inner and 
outer radii are 1 and 10, respectively.  The deformed meshes are generated
by translating the inner circle of the annulus outwards 0.5 units while 
rotating it counterclockwise 10 degrees per timestep.  The mesh is 
deformed until element reversal occurs.  
The meshes resulting from applying various amounts 
of deformation are shown in Figure~\ref{fig:2d_mappings}.  The meshes
in Figure~\ref{fig:2d_mappings}(a) through Figure~\ref{fig:2d_mappings}(c)
are valid meshes.  However, the mesh in Figure~\ref{fig:2d_mappings}(c)
contains elements that are near reversal.  The mesh in Figure~\ref{fig:2d_mappings}(d)
contains reversed elements and is invalid.  The reader is referred to~\cite{ShontzThesis}
for additional examples of triangular and tetrahedral meshes and boundary
deformations which resulted in deformed meshes with reversed
elements.  The purpose of this section is to explore the causes
for reversal and propose some workarounds.  The discussion of
workarounds continues into the next section.

Recall that $\Omega$ is the original polygonal or polyhedral domain, and 
that $\hat\Omega$ is the domain whose boundary is given by the user-specified
deformation of the boundary vertices of $\Omega$, i.e., by the
vertices at coordinates $(\hat x_B,\hat y_B)$.  Assume that this
user-specified deformation is not self-intersecting and preserves
orientation.

Let $\phi$ be the mapping from $\Omega$ to $\hat\Omega$ computed
by FEMWARP.  In other words, interpolate the interior vertex deformations
linearly over the elements to arrive at a continuous function on
the whole domain.  (In the case that FEMWARP fails, parts of
$\phi(\Omega)$ may protrude outside of $\hat\Omega$.)

Let $\phi^*$ be the mapping that is obtained from the exact (continuum)
Laplacian.
In other words, solve the boundary value problems
\begin{eqnarray*}
\bigtriangleup \hat x &=& 0\quad\mbox{on $\Omega$}, \\
\bigtriangleup \hat y &=& 0\quad\mbox{on $\Omega$}, \\
\hat x &= & \hat x_B\quad\mbox{on $\partial\Omega$}, \\
\hat y &= & \hat y_B\quad\mbox{on $\partial\Omega$}
\end{eqnarray*}
and define $\phi^*(x,y)=(\hat x(x,y),\hat y(x,y))$.
Let us call this warping algorithm ``continuum FEMWARP''.

Finally, let $\phi^{+}$ be the mapping that is obtained by
linear interpolation over the elements of $\phi^*$ evaluated
at vertices.  Thus, $\phi^+$ is intermediate between $\phi$ and
$\phi^*$ in the sense
that for $\phi^+$ we use the exact solution to the continuum problem
only at vertex points and use interpolation elsewhere.

There are three possible reasons that FEMWARP could fail:
\begin{enumerate}
\item
Mapping $\phi^*$ might have reversals.
\item
Mapping $\phi^+$ might have reversals even though $\phi^*$
has none.
\item
Mapping $\phi$ might have reversals even though $\phi^+$
has none.
\end{enumerate}
Let us call these Type 1, Type 2, and Type 3 failures.
For Type 1 failure, ``reversal'' means the existence of
a point $x\in\Omega$ such
that $\nabla\phi(x)$ has a nonpositive determinant.  For the second and
third types, ``reversal'' means that a triangle is reversed.  Type 1 
reversals are caused by the boundary deformation alone and are not
related to the mesh. Type 2 reversals are due 
to continuous versus discrete representation of $\phi^*$, 
and Type 3 reversals can be analyzed 
using traditional error estimates for the finite element method.

Let us start with Type 1 reversals.
It is difficult to characterize inputs for which $\phi^*$
will have reversals.  In \cite{ShontzThesis}, a sufficient condition
is given for two-dimensional domains to ensure
that $\phi^*$ will be an invertible function, but the condition 
is unrealistically stringent and is nontrivial to check in practice.

\begin{figure}
\begin{center}
\epsfig{file=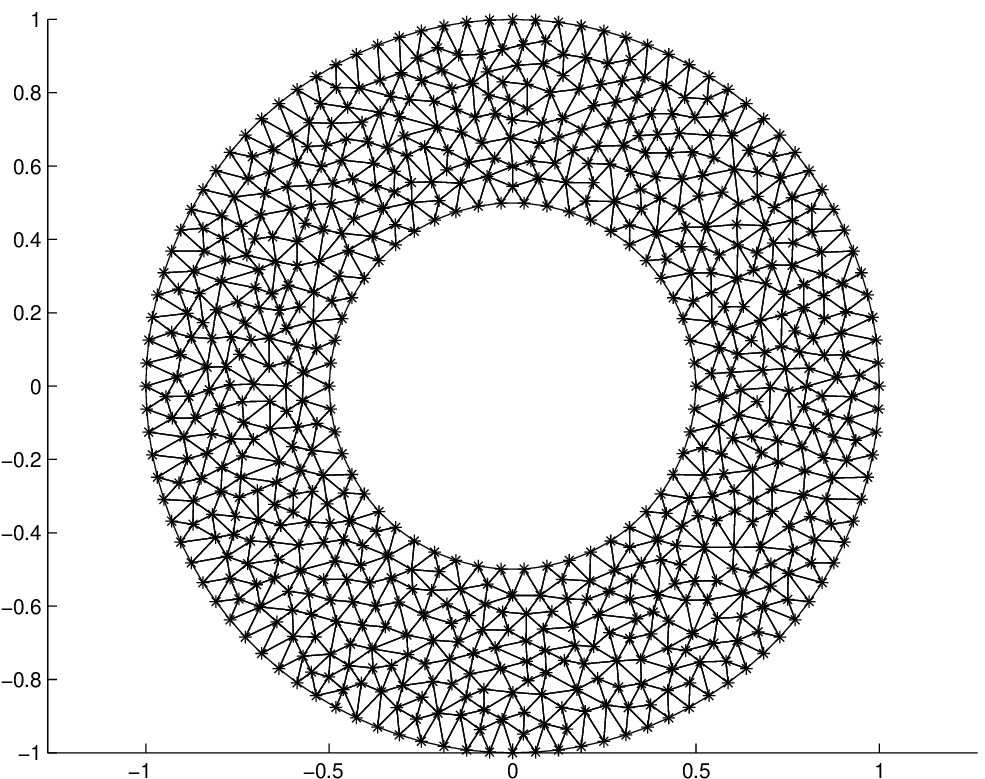,height=4cm}
\end{center}
\caption{All meshes in this section were generated by the
Triangle mesh generator; this is an example of a mesh used
herein.}
\label{fig:annmesh100}
\end{figure}

Rather than presenting the theorem from \cite{ShontzThesis}, we choose 
to present a series of examples
and discussion on how to avoid reversals.
The geometry for most of the examples
in this section is a two-dimensional annulus with outer radius 1 and 
inner radius $r<1$.
Meshes used for tests in this section were generated by 
Triangle, a two-dimensional quality mesh generation package \cite{shew}.
A typical mesh that occurs in 
some of the experiments is depicted in Fig.~\ref{fig:annmesh100};
this mesh has $1238$ elements, $694$ vertices, and 
a maximum side length of $0.1137$.  Its inner radius is $0.5$.

The boundary transformations applied in this section usually
consist of two motions:
First, the inner circular boundary is moved
radially outward (``compression'') to new radius $s$ such that $r\le s<1$.
Second,
a rotation of magnitude $\theta$ 
is applied to the outer boundary.  The relevant
Laplace equations determining $\phi^*$ can be solved in closed form
to yield
$$\phi^*(x,y)=(Ax+By,-Bx+Ay)$$
where $A=a+b/(x^2+y^2)$ and $B=c+d/(x^2+y^2)$, and $a,b,c,d$ are constants
determined by the boundary conditions.  In particular, to
match the boundary conditions just described, one must satisfy
the equations $a+b=\cos\theta$, $c+d=\sin\theta$, $a+b/r^2=s/r$,
$c+d/r^2=0$.  These equations are uniquely solved by choosing
\begin{equation}
\{a,b,c,d\}=\frac{\{\cos(\theta)-rs, rs-r^2\cos(\theta),\sin(\theta),-r^2\sin(\theta)\}}{1-r^2}.
\label{eq:abcd}
\end{equation}

This function $\phi^*$ is invertible, i.e., avoid reversals, provided
the determinant of its Jacobian is always positive.  This determinant
may be computed in closed form:
$\det(\nabla\phi^*)=a^2+c^2-(b^2+d^2)/(x^2+y^2)^2$.  This quantity is minimized when
$x^2+y^2=r^2$.  Therefore, reversals of Type 1 occur if and only if
$r^2(a^2+c^2)\le b^2+d^2$.  Substituting the above formulas
for $a,b,c,d$ yields the result that reversals occur if and only
if
$$2r\cos(\theta)-r^2s-s<0.$$
For example, if $r=s=0.5$ (no compression), then
reversals occur when $\cos(\theta)<.625$, i.e.,
$\theta\ge 51.4...^\circ$.
If $r=.5$ while $s=.75$, then reversals occur
when $\theta \ge 20.4...^\circ$.

We tested the FEMWARP algorithm on the cases described
above with a mesh for the annular region as discussed earlier.
We used a mesh with inner radius $r=0.5$.
This particular mesh contained 10,950 triangles with maximum side length
of 0.039.  For $s=0.5$, when we selected $\theta=51^\circ$, FEMWARP
ran on this mesh without reversals, whereas $\theta=52^\circ$ caused reversals.
As mentioned in the previous paragraph, $\theta\approx 51.4^\circ$ is the
cutoff for Type 1 reversals.  When $r=0.5$, $s=0.75$, FEMWARP
succeeded for $\theta=22^\circ$ but failed for $\theta=23^\circ$, again,
very close to the cutoff for Type 1 reversals.

In addition, we tested the FEMWARP algorithm on two meshes on a cylinder
geometry with height 2 and radius $r = 1$.  The first mesh was a coarse
mesh containing 4320 unstructured tetrahedra which were arranged in
10 layers with 432 tetrahedra per layer.  The second mesh was a fine 
mesh containing 101,550 tetrahedra which were arranged in 25 layers with 
4062 tetrahedra per layer.  Each cylinder mesh was created by first 
using Triangle to generate a mesh on a disk and second using our own
software to extrude the tetrahedral mesh in layers.  The boundary
deformation used for the experiment was $(r \cos(\theta), r
\sin(\theta), z)\mapsto(r \cos(\theta+tz), r \sin(\theta+tz), z),$ where $t$
is a parameter which controls the amount of deformation.  When $t=0$, no deformation 
occurs, whereas increasing values of $t$ correspond to increasing amounts of 
deformation.  In each case, the value of $t$ was increased in increments 
of $0.1$, thus applying an increasing 
amount of deformation to the boundary, until element reversal occurred.
On the coarse mesh, FEMWARP succeeded for $t^{\rm max}_{\rm coarse} = 5.1$
and failed for $t^{\rm fail}_{\rm coarse} = 5.2$, whereas on the fine mesh,
it succeeded for $t^{\rm max}_{\rm fine} = 5$ and failed for 
$t^{\rm fail}_{\rm fine} = 5.1$.  The element
reversal which occurs on the fine mesh for 
$t^{\rm max}_{\rm fine} = 5$ is considered a Type 
1 reversal because the fine mesh is very close to the true Laplace solution, and a Type 
1 reversal means that the true Laplace solution has a Jacobian with a nonpositive 
determinant.  On the other hand, when a reversal occurs on the coarse mesh,
it could be either Type 1, Type 2, or Type 3.  Because 
$t^{\rm max}_{\rm coarse}$ is close to $t^{\rm max}_{\rm fine}$, it is safe to conclude 
that the reversals in the coarse mesh are usually Type 1.

The point of the experiments in the previous two paragraphs is that for
a reasonably refined and reasonably high-quality mesh, most FEMWARP
reversals seem to be Type 1 reversals.
In other words, FEMWARP fails when continuum FEMWARP fails
or is close to failure.
Other experiments not reported here
seem to confirm this point.  Therefore, in order to extend
the range of deformations that can be handled by FEMWARP, the best
strategy is to come up with a way to avoid Type 1 reversals.

One simple method to avoid Type 1 reversals is to take several smaller steps
instead of one big step.  For example, suppose $(\hat x_B',\hat y_B')$ are
positions for the boundary vertices intermediate between their initial
positions and their final positions $(\hat x_B,\hat y_B)$.  Then one could 
define a two-step continuum FEMWARP as follows.  Solve
\begin{eqnarray*}
\bigtriangleup \hat x' &=& 0\quad\mbox{on $\Omega$}, \\
\bigtriangleup \hat y' &=& 0\quad\mbox{on $\Omega$}, \\
\hat x' &= & \hat x_B'\quad\mbox{on $\partial\Omega$}, \\
\hat y' &= & \hat y_B'\quad\mbox{on $\partial\Omega$}
\end{eqnarray*}
for $\hat x'$ and $\hat y'$ to determine a mapping $\phi_1:\Omega\rightarrow
\Omega'$ given by $(x,y)\mapsto(\hat x',\hat y')$ (where $\Omega'$ is
the domain bounded by $(\hat x_B',\hat y_B')$) followed by
\begin{eqnarray*}
\bigtriangleup \hat x &=& 0\quad\mbox{on $\Omega'$}, \\
\bigtriangleup \hat y &=& 0\quad\mbox{on $\Omega'$}, \\
\hat x &= & \hat x_B\quad\mbox{on $\partial\Omega'$}, \\
\hat y &= & \hat y_B\quad\mbox{on $\partial\Omega'$}
\end{eqnarray*}
for $\hat x$, $\hat y$ to obtain a map $\phi_2$.  
Finally, $\phi^*=\phi_2\circ\phi_1$.



The idea in the previous paragraph can be extended to more steps with
smaller increments.  The limiting case of an infinite number of
infinitesimal steps yields an algorithm that we call
``infinitesimal-step continuum FEMWARP.''  Infinitesimal-step FEMWARP
can be described formally as follows.  Let the initial domain be
denoted as $\Omega\subset \R^d$ with boundary $\Gamma$.  Assume there
is a $C^1$ function $g:\Gamma\times\R\rightarrow\R^d$ such that
$g(x,0)\equiv x$; this function $g$ describes the motion of the
boundary.  Here $d=2$ or $d=3$.
Thus, the boundary at time $t$ is denoted as $\Gamma(t)$
and satisfies $\Gamma(t)=g(\Gamma,t)$.  It is assumed that $g(x,t)$ is
injective as a function of $x$, so that the boundary never intersects
itself.

The step at time $t$ is determined by the derivative of $g$.  In particular,
let us define a function $G(x,t)$, where $x\in\Omega(t)$, as follows.
Temporarily fix a particular $t$. 
Solve the vector-valued Laplace equation for $u:\Omega(t)\rightarrow \R^d$
(i.e., a separate Laplace equation
for each coordinate entry) given by 
$$
\begin{array}{rcll}
\bigtriangleup u(x) &=& 0 & \mbox{for $x\in \Omega(t)$,} \\
u(x) & = & \frac{\partial}{\partial t} g(g^{-1}(x,t),t) & \mbox{for $x \in \Gamma(t)$}.
\end{array}
$$
Finally, define $G(x,t)=u(x)$.  Last, given a point $x_0\in \Omega$, we consider
the trajectory defined by the initial value problem $x(t)=x_0$; $x'(t)=G(x,t)$.
The solution operator of this ODE system, say $\Phi(x,t)$, defines the
infinitesimal-step continuum FEMWARP mapping at time $t$.  

It can be shown that this map is a bijection with a positive Jacobian for all $x$
(i.e., no reversals).  This is a consequence of the well-known standard fact
that the solution operator for an ODE system with a Lipschitz forcing function $G$
is bijective with a positive Jacobian.  The usual textbook theorem does not quite apply
to this case because the spatial domain of $G(x,t)$ is not fixed in time but depends
on $t$.  However, the theorem is still valid because for any $x_0$ interior to $\Omega(t)$,
one can break up the trajectory into small pieces and define $G(x,t)$ only in a small
neighborhood around $x$ (but fixed over time for each piece).  Then the pieces can be
assembled together to prove the result.

In the case of the annulus, it is possible again to write down
infinitesimal-step continuum FEMWARP in closed form.  For simplicity,
let us assume $r=s$ so that the only deformation is the rotation of the
outer boundary.  Assume this rotation is broken up into infinitesimally
small rotations.  (Another choice would be to connect the initial positions
to the final positions with line segments, and break up the boundary motion as
infinitesimal increments along the line segments.  This way to obtain
a continuous boundary 
motion is undesirable, however, because for a sufficiently
large rotation, the line segments would cut through the inner boundary 
of the annulus and hence
cause tangling of the boundaries.)

With the setup described in the last paragraph, the deformation
for an outer rotation of $\theta$ 
computed by infinitesimal-step continuum FEMWARP maps a point at
initial position $\rho(\cos\phi,\sin\phi)$ ($r\le\rho\le 1$) to
$\rho(\cos(\phi+\alpha),\sin(\phi+\alpha))$, where
$\alpha=(1-r^2/\rho^2)\theta/(1-r^2)$.  This map is clearly bijective
for any value of $\theta$; it corresponds to rotating each concentric
circle of the annulus by an amount that interpolates between 0 (when
$\rho=r$) and $\theta$ (when $\rho=1$).

Thus, by using small-step FEMWARP with sufficiently small steps,
we can essentially eliminate Type 1 failures. Small-step FEMWARP preserves
the attractive property of FEMWARP that it is exact for affine maps,
as long as all the intermediate steps are also affine.  Unfortunately,
it loses the
attractive property that only one coefficient
matrix for solving the linear system needs to be
factored.  Small-step FEMWARP requires the solution of a different 
coefficient
matrix for each step.  This drawback is partly ameliorated by the fact
that even though the matrices are different, they have the same nonzero
pattern, and hence the symbolic phase of sparse direct solution may be
reused.  If instead an iterative method is being used to solve the mesh 
warping equations, then the sparsity pattern may be reused in the 
preconditioner.  In addition, the factored coefficient matrix at step 
$t_k$ can be used as a preconditioner for an iterative method at step
$t_{k+1}$.

Elimination of Type 1 failures means that 
the mapping function $\phi^*$ has no reversals in the sense that
the determinant of its Jacobian is positive everywhere; equivalently,
it does
not reverse any infinitesimally small triangles.  
A Type 2 failure occurs because the triangles in the mesh
have finite (non-infinitesimal) size and hence can still be
reversed by $\phi^+$. The following 
theorem
characterizes when this can happen.

\begin{theorem}
Suppose that 
$f: \Omega \rightarrow \R^2$ is bijective, orientation-preserving 
and $C^2$ on $\Omega$
with $\nabla f$ nonsingular.  
Let $T$ be a triangle in the mesh with vertices $\{v_1,v_2,v_3\}$,
and let $T'$ be the triangle whose vertices are $\{f(v_1),f(v_2),f(v_3)\}$.
If \eref{eq:revineq} below holds,
then $T'$ is not reversed.
\label{thref:inv}
\end{theorem}
\begin{proof}
Recall that  triangle $T$ with vertices $\{v_1,v_2,v_3\}$ is positively
oriented if and only if $\det(A)>0$, where
$$A=(v_2-v_1,v_3-v_1).$$
In order to analyze the analogous quantity 
for $\{f(v_1),f(v_2),f(v_3)\}$,
we start with the following algebra, which invokes the fundamental
theorem of calculus twice:
\begin{eqnarray*}
f(v_2)-f(v_1) 
& = & 
\int_0^1 \nabla f((1-t)v_1+tv_2)(v_2-v_1)\,dt \\
& = & 
\left(\int_0^1 \nabla f((1-t)v_1+tv_2)\,dt\right)(v_2-v_1) \\
& = & 
\left( \nabla f(v_1)+
\int_0^1 \left[\nabla f((1-t)v_1+tv_2)-\nabla f(v_1)\right]\,dt\right)
(v_2-v_1) \\
& = & 
\nabla f(v_1)(v_2-v_1)\\
& & +\left(\int_0^1 \left[ \int_0^1 \nabla^2 tf((1-s)v_1+s((1-t)v_1+tv_2))(v_2-v_1)\,ds\right]\,dt \right) \\
& & \ \ \ \ \ \ \ \ \ \ \ \ \ (v_2-v_1) \\
& = & 
\nabla f(v_1)(v_2-v_1)+e_1
\end{eqnarray*}
where $\Vert e_1\Vert \le Mh^2$, where $h$ is the maximum side length of
$T$ (an upper bound on $\Vert v_2-v_1\Vert$) and $M$ is an upper bound
on $\Vert \nabla^2 f\Vert$ in the triangle.  Similarly,
$$f(v_3)-f(v_1)=\nabla f(v_1)(v_3-v_1)+e_2,$$ 
where again $\Vert e_2\Vert \le Mh^2$.  
Therefore,
\begin{equation}
(
f(v_2)-f(v_1),
f(v_3)-f(v_1)
) = \nabla f(v_1)A + E
\label{eq:pertAE}
\end{equation}
where $A$ is as above and $\Vert E\Vert_2\le \sqrt{2}Mh^2$.
Observe that $\nabla f(v_1)A$ has positive determinant by assumption.
Therefore, the left-hand side can have negative determinant only if
$E$ is a sufficiently large perturbation to change the determinant
sign.  If $E$ is such a large perturbation, then by the continuity
of the determinant, there is a perturbation
$E'$ no larger than $E$ such that the $\det(\nabla f(v_1)A+E')=0$,
i.e., $\nabla f(v_1)A+E'$ is singular.  Furthermore, $\Vert E'\Vert
\le \sqrt{2}Mh^2$.  By Theorem 2.5.3 of \cite{GVL3}, this means that
\begin{equation}
\sqrt{2}Mh^2\ge \sigma_{\min}(\nabla f(v_1)A)\ge \sigma_{\min}(\nabla f(v_1))
\sigma_{\min}(A),
\label{eq:Mh2}
\end{equation}
where $\sigma_{\min}(A)$ and $\sigma_{\max}(A)$ denote the smallest 
and largest singular values of $A$, respectively.
It follows from the equation $AA^{-1}=I$ that 
that the columns of $A^{-1}$ are parallel to the  altitude segments of triangle
$T$ perpendicular to $v_1v_3$ and $v_1v_2$ respectively, but scaled
so that their lengths are the reciprocals of those altitude lengths.
Therefore, $\sigma_{\max}(A^{-1})\le \sqrt{2}/\minalt(T)$, where
$\minalt(T)$ means the minimum altitude.  Thus, $\sigma_{\min}(A)\ge
\minalt(T)/\sqrt{2}$.  Substituting this inequality into
\eref{eq:Mh2} and rearranging yields
$$\frac{\sigma_{\min}(\nabla f(v_1))}{M}\le 2h\asp(T)$$
where $\asp(T)$, the {\em aspect ratio} of $T$, equals $h/\minalt(T)$.  The
aspect ratio is often used as a shape-quality metric; lower values mean
a better shaped triangle.  Thus, reversal cannot happen if the
opposite inequality holds:
\begin{equation}
\frac{\sigma_{\min}(\nabla f(v_1))}{M}> 2h\asp(T).
\label{eq:revineq}
\end{equation}
\end{proof}

The point of this theorem is that Type 2 reversals cannot occur
for a sufficiently refined mesh (i.e., $h$ sufficiently small in
\eref{eq:revineq}), assuming that the mesh quality does not
decay, and assuming that $\phi^*$ is a nonsingular function.
(Assuming Type 1 reversals are excluded, 
function $\phi^*$ is never singular on the interior because Laplace
solutions are analytic.  It could be singular at the boundary if,
for example, $\Omega'$ has a corner where $\Omega$ had none.)

We tested this theorem for two examples, each of which diverges a bit from
the theoretical prediction. 
For the first example, we generated a uniform mesh for the rectangle
$[0,2]\times[0,1]$ using Triangle and mapped all the vertices 
using the function $f(x,y)=(x,y+\alpha x(2-x))$.  For each
mesh, $\alpha$ was incremented by 1 until reversal occurred.   (No Laplace
solution was involved in this test case.)  We tabulated the
values of $h$ versus $\alpha$ in Table~\ref{tab:rectable}.  As predicted
by the theorem, the table shows that as $h$ decreases, a larger value
of $\alpha$ is tolerated.  Contrary to the theorem, however, the table
shows that $\sigma_{\min}(\nabla f)/\Vert\nabla^2 f\Vert$ is decreasing
faster than $h$.  In other words, reversals are avoided to a greater
extent than predicted by the theorem.
The reason  for this discrepancy is that the
perturbation term $E$ in
\eref{eq:pertAE} is not well aligned with the direction that
drives $(\nabla f)A$ toward singularity in this example.  In particular, $E$
affects only the $y$-components (since the transformation is
linear for $x$-coordinates).  On the other hand, transformation 
$\phi$ stretches the
triangles substantially in the $y$-direction, so that the most effective
way to perturb $(\nabla f)A$ toward singularity
is a small change to the $x$-components.

As a second test case, consider the transformation of the annulus with
radii $(0.5, 1)$
that results from continuous-warping continuum FEMWARP, that is, the
transformation that rotates a point at radius $\rho$ by
angle $\alpha(1-r^2/\rho^2)/(1-r^2)$, where $r$ is the inner radius
($r=0.5$ for this test).  For each mesh, the parameter $\alpha$ was stepped
in increments of $\pi/16$ until reversals were encountered.
We tabulated values of $h$ versus the first
of $\alpha$ causing failure in Table~\ref{tab:anntable}. 
As predicted by the theorem, decreasing $h$ corresponded to increasing
values of $\alpha$, i.e., greater distortion of the domain.
Again, these results do not initially
 correspond to the preceding theorem quantitatively:
in this case, $h$ is decreasing faster than
$\sigma_{\min}(\nabla f)/\Vert\nabla^2 f\Vert$.  Only the last three
rows of the table show that $h$ and 
$\sigma_{\min}(\nabla f)/\Vert\nabla^2f\Vert$ are decreasing proportionally.
The reason for this discrepancy is that $\nabla^2 f$ is much larger
on the inner boundary than elsewhere, and the meshes in the initial
rows of the table are not sufficiently refined to resolve the 
variation in the value of the derivative.

\begin{table}
\begin{center}
\caption{Second column $\alpha^{\rm fail}$ is the first value
of $\alpha$ in the transformation $(x,y)\mapsto(x,y+\alpha x(2-x))$
of a rectangle
that causes reversals in the mesh.  The first column shows
the mesh cell size (maximum edge length). The third column
is  $\sigma_{\min}(\nabla f)/\Vert \nabla^2 f\Vert$ evaluated at
a vertex of a triangle that reversed.}
\label{tab:rectable}
\begin{tabular}{ccc}
\hline
$h$ & $\alpha^{\rm fail}$ & $\sigma_{\min}(\nabla f)/\Vert \nabla^2 f\Vert$ \\
\hline
0.205 & 10 & $4.6\cdot 10^{-3}$\\
0.108 & 14 & $2.5 \cdot 10^{-3}$\\
0.057 & 30 & $5.7 \cdot 10^{-4}$\\
0.030 & 51 & $1.9 \cdot 10^{-4}$\\
0.015 & 95 & $5.6\cdot 10^{-5}$\\
\hline
\end{tabular}
\end{center}
\end{table}

\begin{table}
\begin{center}
\caption{Second column $\alpha^{\rm fail}$ is the first value
of $\alpha$ in the transformation 
$(\rho,\theta)\mapsto(\rho,\theta+\alpha(1-r^2/\rho^2)/(1-r^2))$ (in
polar coordinates) of an annulus
that causes reversals in the mesh.  The first column shows
the mesh cell size (maximum edge length). The third column
is  $\sigma_{\min}(\nabla f)/\Vert \nabla^2 f\Vert$ evaluated at
a vertex of a triangle that reversed.}
\label{tab:anntable}
\begin{tabular}{ccc}
\hline
$h$ & $\alpha^{\rm fail}$ & $\sigma_{\min}(\nabla f)/\Vert \nabla^2 f\Vert$\\
\hline
0.202 & $7\pi/16$ & $6.9\cdot 10^{-3}$\\
0.114 & $\pi/2$ & $4.6\cdot 10^{-3}$\\
0.058 & $9\pi/16$ & $6.2\cdot 10^{-3}$\\
0.031 & $5\pi/8$ & $2.5\cdot 10^{-3}$ \\
0.015 & $11\pi/16$ & $2.6\cdot 10^{-3}$ \\
0.008 & $13\pi/16$ & $1.4\cdot 10^{-3}$ \\
0.004 & $\pi$ & $7.2\cdot 10^{-4}$ \\
\hline
\end{tabular}
\end{center}
\end{table}

Thus, we have seen that Type 1 reversals can be avoided by using
small-step FEMWARP instead of FEMWARP, and Type 2 reversals can
be avoided by using a sufficiently refined mesh.  The remaining type
of reversals, Type 3, are rare according to our experiments.  
Type 3 reversals are caused by the difference between the true
value of the Laplace solution and the finite element approximation
to that solution.  Intuitively, this phenomenon should not be commonplace
because the perturbation size to a mesh necessary to cause a reversal
of a triangle of side-length $h$ is $O(h)$, whereas the difference between
the two mappings is $O(h^2)$.

Consider the following test.  We generated a sequence of small-step
rotations
for the annulus in two ways.  In the first case, we took steps that rotate
the outer boundary by $\pi/16$ and leave the inner boundary invariant, 
each time computing the Laplace solution exactly analytically.
This corresponds to iteratively applying the transformation 
$\phi^*(x,y)=(Ax+By,-Bx+Ay)$ to the mesh, 
where $A,B,C,D$ are functions of $x^2+y^2$ as above, and $a,b,c,d$ are constants 
determined by \eref{eq:abcd} with $r=s=0.5$, $\theta=\pi/16$.
In the second case, we solved
Laplace's equation for the above boundary condition using the finite element
mesh that results from small-step FEMWARP.  In both cases we tried meshes
with several different values of $h$.  The results are tabulated in
Table~\ref{tab:typeiiitest}.  As can be seen, discretized
small-step FEMWARP outperformed continuum small-step FEMWARP.  

In other
words, not only did Type 3 reversals not occur, but in fact it
seems to be preferable to use the discretized solution for mesh warping
rather than the continuum solution.
The difference between $\phi^+$ and $\phi$ is the
usual discretization error in finite element methods.   A possible
explanation for the improved resistance to reversals of the
finite element solution is as follows.  After several steps of 
small-step FEMWARP, Laplace's equation is  solved on a mesh with 
mostly poorly-shaped elements, some extremely poorly-shaped. 
A Laplace solution minimizes the functional $F(u)=\int_\Omega \nabla u \cdot
\nabla u$ over $H^1$ functions $u$ on 
the domain, and the finite element solution minimizes the
same functional $F(u)$ over the space of piecewise linear choices for 
$u$ \cite{johnsonbook}.  
A very
poorly-shaped element is ``stiffer'' than others in the following sense.
An affine linear
function $u$ defined over a triangle that has an angle close to $180^\circ$
will have a quite large gradient value (compared to a well-shaped triangle
with the equal area and equal vertex values) 
unless the vertex values lie in a certain restricted 
range.  Therefore, the extra stiffness of these elements will cause them
to be deformed less than the better-shaped elements in the
optimal solution that minimizes $F$.  Since the poorly-shaped
elements are those most in danger of 
being reversed, this is a desirable effect.

\begin{table}
\begin{center}
\caption{Continuum versus discretized small-step FEMWARP applied to
a mesh of an annulus in order to test for Type 3 reversals.
The table shows that discretized small-step FEMWARP seems less
prone to reversals than continuum small-step FEMWARP.}
\label{tab:typeiiitest}
\begin{tabular}{ccc}
\hline
$h$ & $\alpha^{\rm fail}_{\rm continuum}$ & $\alpha^{\rm fail}_{\rm discrete}$ \\
\hline
$0.202$ & $7\pi/16$ &  $\pi/2$ \\
$0.114$ & $\pi/2$ & $9\pi/16$ \\
$0.058$ & $9\pi/16$ & $5\pi/8$ \\
$0.031$ & $5\pi/8$ & $7\pi/8$ \\
$0.015$ & $11\pi/16$ & $\pi$ \\
\hline
\end{tabular}
\end{center}
\end{table}

The next topic to consider in this section is how to select a stepsize
for small-step FEMWARP.  The theory developed above indicates that as long
as the step size is well below a step large enough to cause reversals
of $\phi^*$, the step size should not matter so much.  In fact, we propose
the following simple strategy, which seems
to be effective.  Attempt to take a 
very large step (e.g., a rotation of size $\pi$ in the case of the annulus).
If this fails (causes reversals), then halve the stepsize and try again
until success.  Update the mesh and try another such step.  Note that
in the process of searching for a correct stepsize, the coefficient
matrix in FEMWARP is the same for each trial.  Therefore, the Cholesky
factors do not need to be recomputed until the mesh is updated.

Another way to carry out small-step FEMWARP would be to take 
constant (small) steps on each iteration.  We compared these two methods and
found that the first was much more efficient, and furthermore, 
reversals are resisted better by the first
strategy.  Therefore, the repeated halving strategy
is recommended.  Table~\ref{tab:stepsizetest} summarizes the result
for the annulus again.  For the halving strategy, updates were pursued
until the stepsize dropped below $\pi/128$.  For the constant-step strategy,
the stepsize was taken to be $\pi/128$.

\begin{table}
\begin{center}
\caption{Variable stepsize versus constant stepsize for
small-step FEMWARP applied to a mesh of an annulus.
Second and third columns are the maximum amount of
rotation prior to reversals for the two methods.
Fourth and fifth columns are the number of Cholesky
factorizations required by the two methods.}
\label{tab:stepsizetest}
\begin{tabular}{ccccc}
\hline
$h$ & $\alpha^{\rm max}_{\rm VS}$ & $\alpha^{\rm max}_{\rm CS}$ 
& $\mbox{NCHOL}_{\rm VS}$ & $\mbox{NCHOL}_{\rm CS}$ \\
\\
\hline
$0.202$ & $1.7426$ &  $1.7671$ & $13$ & $72$  \\
$0.114$ & $2.2089$ & $2.0862$ & $24$ &  $85$ \\
$0.058$ & $2.6998$ & $2.4789$ & $29$ & $101$ \\
$0.031$ & $3.4852$ & $2.7980$ & $34$ & $114$ \\
\hline
\end{tabular}
\end{center}
\end{table}

The final workaround to element reversal we consider in this section is 
to use, instead of the discrete harmonic map, the mean value 
map~\cite{floater2} for computing the weights.  The mean value map 
satisfies the same affine exactness as the discrete harmonic map (as shown 
in Lemma 3.1), and because of this may tend to give similar results 
to the discrete harmonic map.  Yet, at the same time, the mean value map 
uses weights that are always nonnegative (regardless of the shapes
of the triangles), and so at least for a convex boundary, it guarantees an 
injective map~\cite{floater2}.

In order to compare the performance of FEMWARP to that of the algorithm 
using the mean value map, two experiments were performed.  The first 
experiment was on an annulus mesh with radii $(0.5,1)$ which was 
composed of 10840 triangles.  A rotation of magnitude $\theta$ was 
applied to the outer boundary.  For both algorithms, $\theta = 50^{o}$ was 
the maximum amount of deformation which could be successfully applied; a 
deformation with $\theta = 51^{o}$ yielded element reversal.  Thus, the 
two algorithms did not realize a difference in performance for this experiment.

The second experiment was on an annulus mesh with radii (0.3,1) 
which was composed of 13162 triangles.  For this experiment, the boundary 
deformation that rotates a point at radius $\rho$ by angle 
$\alpha(\rho-r)/(1-r)$, where $r$ is the inner radius ($r = 0.3$ 
for this test).  This corresponds to rotating each concentric circle of 
the annulus by an amount that interpolates between 0 (when $\rho = r$) and 
$\alpha$ (when $\rho = 1$). The parameter $\alpha$ was stepped in 
increments of $0.05$ until reversals were encountered.  In this case, 
FEMWARP was able to perform 52 rotation steps successfully, whereas the 
mean value map algorithm only performed 24 steps before element reversal 
occurred.  Thus, the discrete harmonic map used in FEMWARP generates 
weights which allow for successful application of larger deformations than 
those tolerated by the mean value map.

\section{Mesh warping and mesh untangling}
\label{sec:untangling}

In the previous section we considered reasons why FEMWARP can fail
and also some possible workarounds.  Some of the workarounds in the 
previous section, however, are not available in all circumstances.  For 
example, the workaround for Type 1 reversals, namely, small-step FEMWARP, 
requires a homotopy from the old to new boundary conditions and also
requires solution of many linear systems with distinct coefficient
matrices.  Although linear interpolation could be performed if a
homotopy from the old to new boundary conditions is not available,
FEMWARP probably would not give the desired motion, e.g., if the 
boundary deformation involves some kind of rotation.
The workaround for Type 2 reversals requires refined meshes, which may 
not be available.

Another workaround is to switch to a different algorithm, for example,
Opt-MS, a mesh untangling method due to Freitag and 
Plassmann~\cite{FreitagPlassmann}.   Opt-MS takes
as input an arbitrary tangled mesh and a specification of which vertices
are fixed (i.e., boundary vertices) and which are movable.  It then
attempts to untangle the mesh with a sequence of individual vertex
moves based on linear programming.  More details are provided below.

In this section, we consider the use of Opt-MS for mesh warping.  We find
that the best method is a hybrid of FEMWARP and Opt-MS.

To untangle the  mesh, Opt-MS performs repeated sweeps over
the interior vertices.  For each interior vertex, it repositions the vertex 
at the coordinates that maximize the minimum signed area
(volume) of the elements adjacent to that vertex
(called the ``local submesh'').  The signed area is negative for a
reversed element, so maximizing its minimum value is an attempt
to fix all reversed elements.

Let $\x$ be the location of the {\em free vertex}, that is,
the current interior vertex being processed in a sweep.
Let $\x_1, \ldots, \x_p$ be
the positions of its adjacent vertices, and $t_1, \ldots, t_n$ be the
incident triangles (tetrahedra) that compose
the local submesh.  Then the function
that prescribes the minimum area (volume) of an element in the local
submesh is given by $q(\x) =\min_{1 \leq i \leq n} \ A_i(\x)$, where $A_i$
is the area (volume) of simplex $t_i$.  In 2D, the area of triangle $t_i$
can be stated as a function of the Jacobian of the element as follows:
$A_i = \frac{1}{2} \det (\x_i-\x, \ \x_j-\x)$ which is a linear 
function of the free vertex position; the same is true in 3D.  Freitag 
and Plassmann use this fact to formulate the solution to $\max 
q(\x) = \max \min_{1 \leq i \leq n} \ A_i(\x)$ as a linear programming 
problem which they solve via the simplex method.  On each sweep, $m$ 
linear programs are solved which sequentially reposition each interior 
vertex in the mesh.  Sweeps are performed until the mesh is untangled or a 
maximum number of sweeps has occurred.  

A shortcoming of Opt-MS, in comparison to FEMWARP, is that it is
not intended to handle a very large boundary motion even
if that motion is affine.  For example, starting from a 2D mesh, if
the boundary vertices are all mapped according to the function 
$(x,y)\mapsto(-x,-y)$ while the interior vertices are left unmoved,
in many cases Opt-MS is unable to converge.   
FEMWARP, on the other hand, will 
clearly succeed according to Theorem~\ref{thm:affine}.

Therefore, we propose the following algorithm: one applies first
the FEMWARP mesh warping algorithm, and if the mesh is tangled,
one uses Opt-MS to untangle the mesh output by
FEMWARP (as opposed to the original mesh).  The rationale for this
algorithm is that FEMWARP is better able to handle the gross
motions while Opt-MS handles the detailed motion better.

To test whether this algorithm works, we applied it again to the mesh
depicted in Fig.~\ref{fig:annmesh100}.  The boundary motion is as
follows: we rotate the outer circular boundary by $\theta_1$ degrees
and the inner boundary by $\theta_2$ degrees and then test three
algorithms: FEMWARP alone, Opt-MS alone, and hybrid.  (The hybrid was
tested only in the case that the two algorithm individually both
failed.)  The results are tabulated in Table~\ref{tab:hybrid_test}.
The table makes it clear that the hybrid often works when Opt-MS and
FEMWARP both fail.  Thus, the hybrid method is another technique for
situations when FEMWARP alone fails and small-step FEMWARP combined
with mesh refinement may be unavailable.  
The hybrid method has the disadvantage, when compared to FEMWARP,
that it does not produce a
continuous motion of interior vertices but rather only a final
configuration.

\begin{table}
\begin{center}
\caption{The row header indicates $\theta_1$ (degrees), the
rotation of the outer boundary,
and the column header $\theta_2$ is the rotation
of the inner boundary.  The table entries are
as follows: `F' means FEMWARP succeeded but not Opt-MS, `O' means
Opt-MS succeed but not FEMWARP,`B' means both succeeded,
and `H' means neither succeeded, but the hybrid succeeded,
and finally `--' means none succeeded.}
\label{tab:hybrid_test}
\begin{tabular}{l|lllllllllllll}
\hline
$\theta_1\backslash\theta_2$ &
$0$ & $15$& $30$ & $45$ & $60$ & $75$ & $90$ & $105$ & $120$ & $135$ & $150$
& $165$ & $180$ \\
\hline
$0$  & B & B & B & B & O & O & O & --& --& --& --& --& -- \\
$15$ & B & B & B & B & B & O & O & H & --& --& --& --& -- \\
$30$ & B & B & B & B & B & B & H & O & H & --& --& --& -- \\
$45$ & B & B & B & B & B & B & B & O & O & H & --& --& -- \\
$60$ & O & B & B & B & B & B & B & B & O & O & H & --& -- \\
$75$ & O & O & B & B & B & B & B & B & B & H & H & H & -- \\
$90$ & O & O & O & B & B & B & B & B & B & B & O & H & H  \\
$105$& O & O & O & O & B & B & B & B & B & F & B & O & H  \\
$120$& --& --& H & H & H & B & B & B & F & F & B & F & H  \\
$135$& --& --& --& H & H & H & F & B & B & B & B & B & B  \\
$150$& --& --& --& --& H & H & H & F & B & B & B & B & B  \\
$165$& --& --& --& --& --& H & H & H & F & F & B & B & B  \\
$180$& --& --& --& --& --& --& H & H & H & F & B & B & B \\
\hline
\end{tabular}
\end{center}
\end{table}

\section{Three-dimensional tests}
\label{sec:test3d}
\begin{figure}[ht]
\begin{center}
\caption{The 3D test meshes used in this section:  (a) foam5~\cite{pat_mesh}, (b) gear~\cite{pat_mesh}, 
(c) hook~\cite{pat_mesh}, and (d) tire~\cite{lori_mesh}.}
\label{fig:test_meshes_3D}
\vbox{
\hbox{
\hspace{.5in}
  \subfigure[foam5]{
    \includegraphics[height=5cm,width=4cm,trim=5 5 5 5, clip=true]{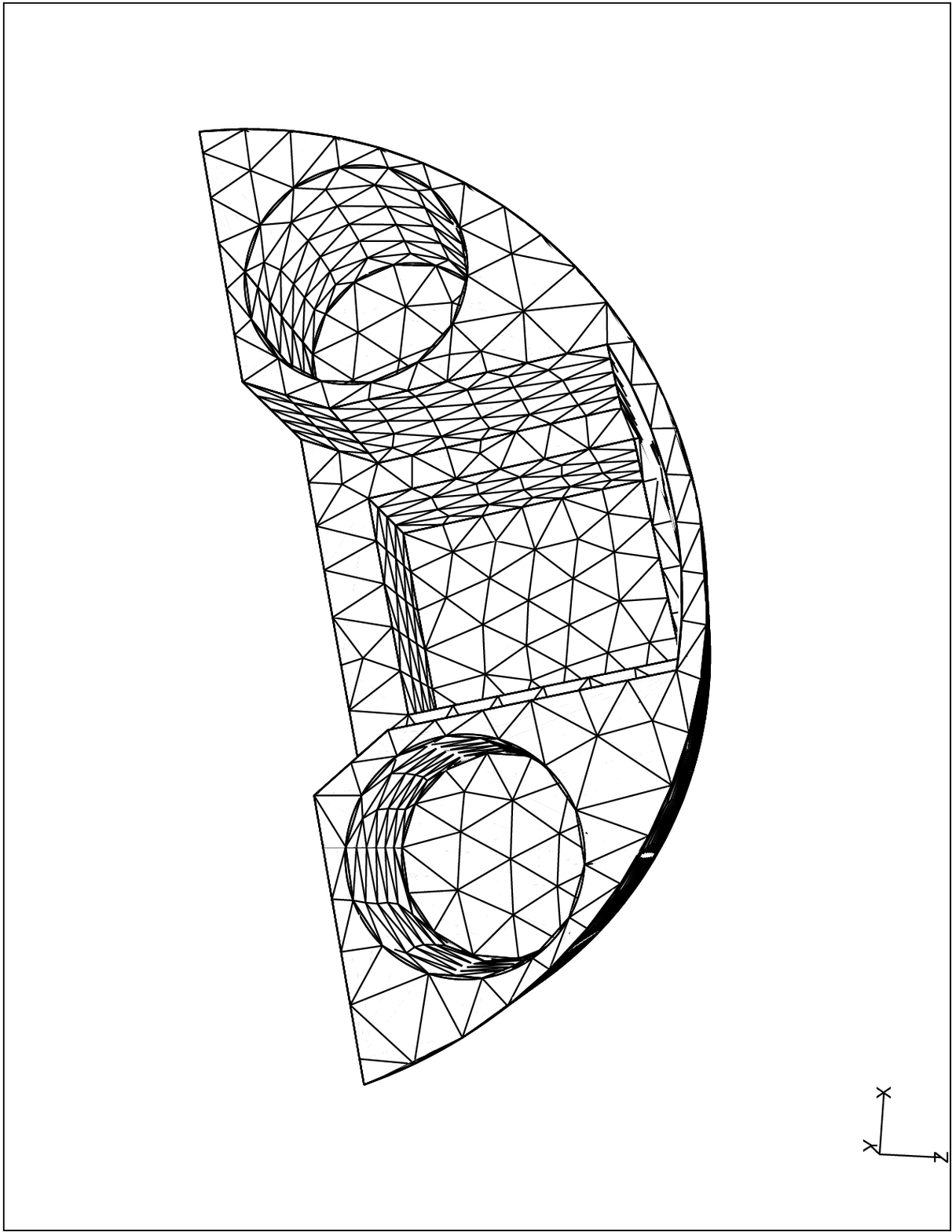}
    \label{fig:foam5}
  }
  \hskip 1cm
  \subfigure[gear]{
    \includegraphics[height=5cm,width=4cm, trim=5 5 5 5, clip=true]{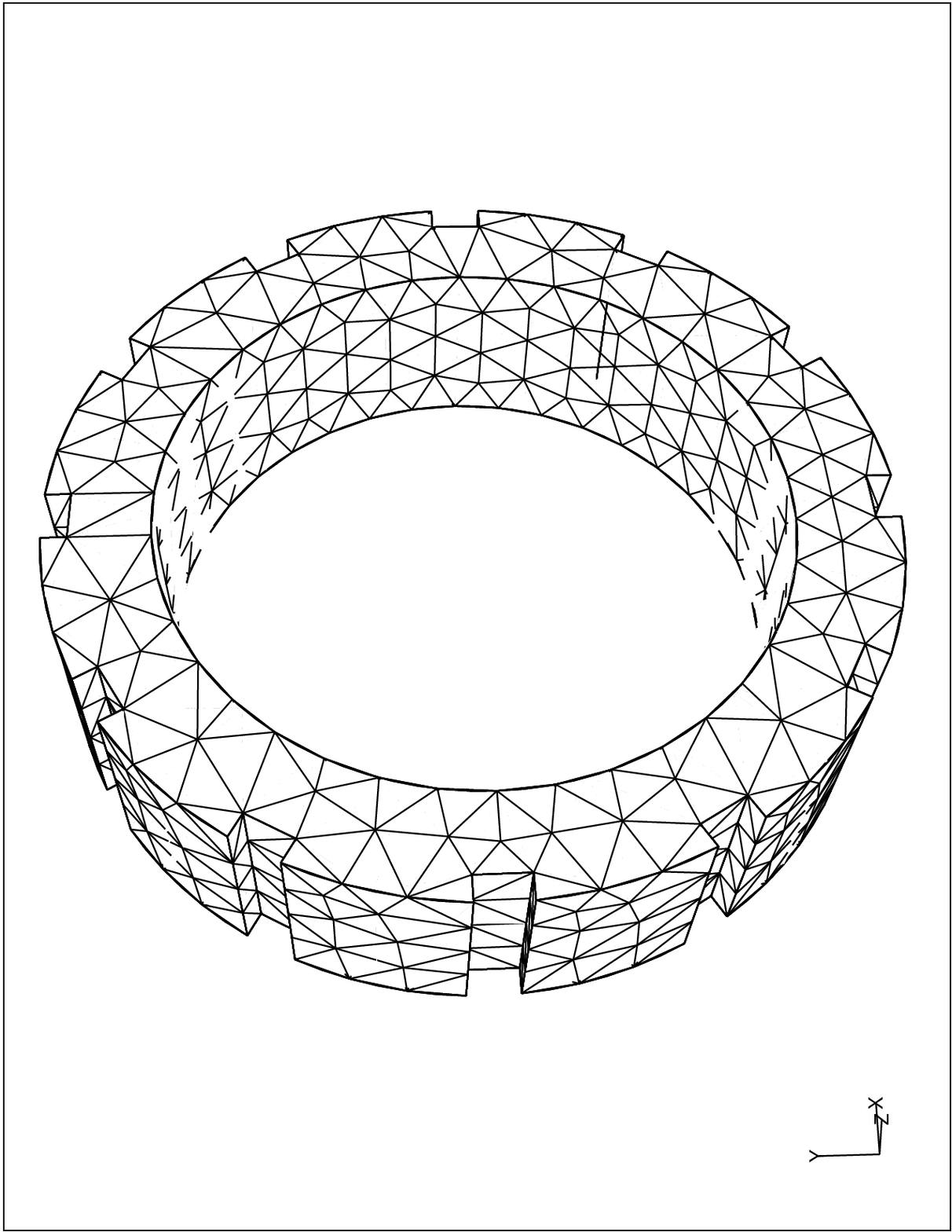}
    \label{fig:gear}
  }
}
\hbox{
\hspace{.5in}
\subfigure[hook]{
\includegraphics[height=5cm,width=4cm,trim= 5 5 5 5, clip=true]{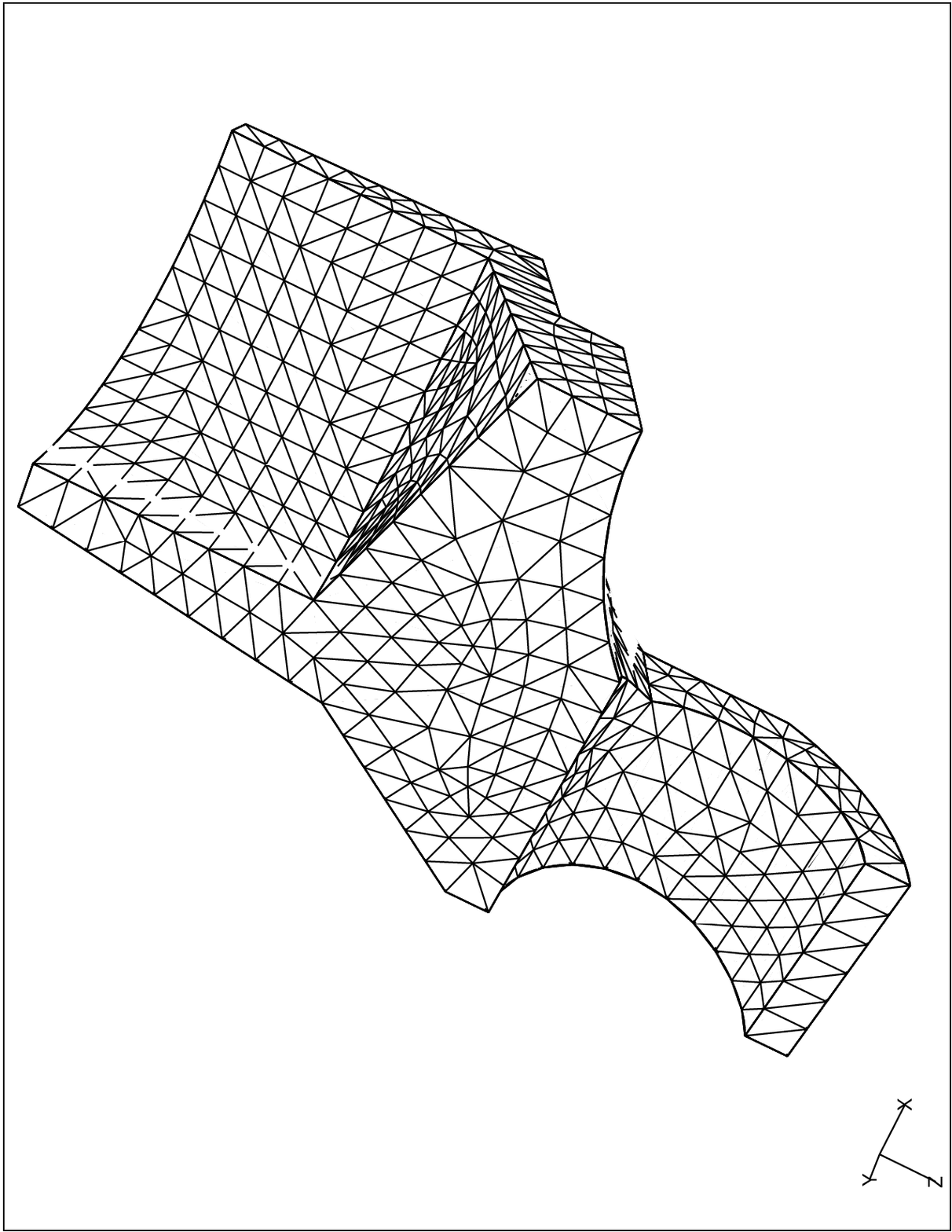}
\label{fig:hook}
}
\subfigure[tire]{
\includegraphics[height=5cm,width=5cm]{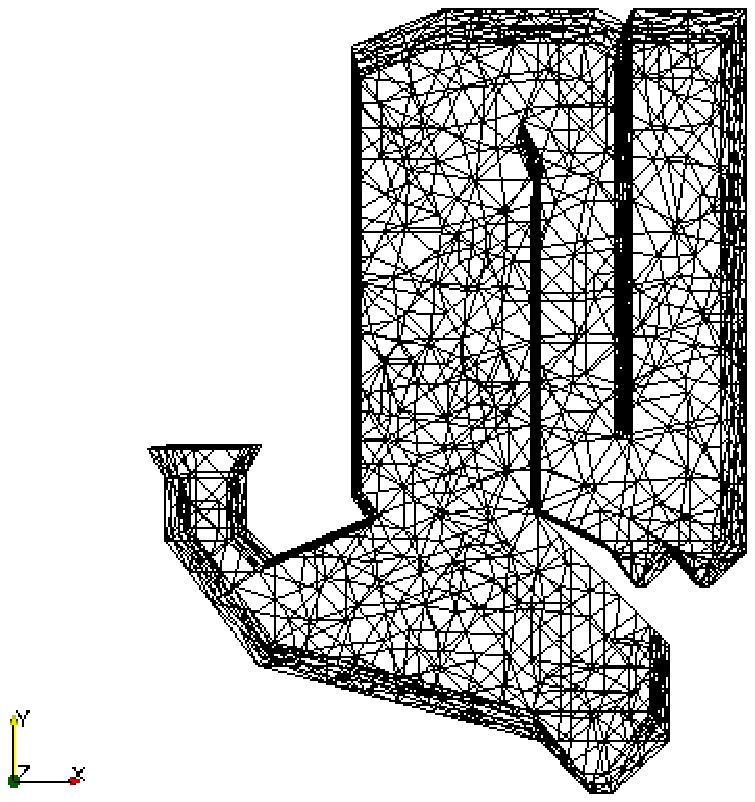}
\label{fig:tire}
}
}
}
\end{center}
\end{figure}

In this section we compare the robustness 
against reversals of FEMWARP, small-step
FEMWARP, and the hybrid FEMWARP/Opt-MS method
on examples of 3D meshes~\cite{pat_mesh} and~\cite{lori_mesh}
shown in Figure~\ref{fig:test_meshes_3D}.

We choose specific 
nonlinear boundary deformations parameterized
by a scalar $\alpha$ in order to determine how much 
deformation each test mesh could withstand when warped according to each 
method:  
$$
\left( \begin{array}{c} x \\ y  \\ z \end{array} \right) \mapsto
\left( \begin{array}{rrr} 2 & -1 & 0 \\ -2 & 5 & 0 \\ 0 & 0 & 1 \\ 
\end{array} \right)
\left( \begin{array}{c} x \\ y \\ z \end{array} \right)
+ \alpha \left( \begin{array}{c} 0.1 x y  \\ 0.5 yz \\ 0.1x^2 \end{array} 
\right).$$

\tabref{tab:compare_3D_1} gives the results obtained 
from warping various three-dimensional 
meshes according to this
boundary deformation.  The data in this table is as follows.
The first three columns give the name of the mesh, the number of
boundary vertices, and the number of total vertices.  The fourth column 
$\alpha_{\rm FEMWARP}^{\rm max}$ gives the maximum
 value of $\alpha$ encountered for which FEMWARP
succeeded.  Parameter $\alpha$ was stepped by $\Delta \alpha$,
where $\Delta\alpha$ is indicated in the last column of the table.
\tabref{tab:compare_3D_2} is a continuation of the previous table.
The second and third columns of this table
show $\alpha_{\rm ssFEMWARP}^{\rm max}$, which is
the last value of $\alpha$ for which small-step FEMWARP succeeds,
and $\mbox{NCHOL}_{\rm SS}$, which is the number of Cholesky factorizations
required by small-step FEMWARP. 
The variable-step version of FEMWARP described in Section~\ref{sec:inversion}
was used.
The minimum allowed stepsize for small-step FEMWARP was also
taken to be $\Delta \alpha$.

The table indicates that small-step FEMWARP was about equal
in robustness
against reversals compared with FEMWARP except
for the last mesh.  In the case of ``tire'', small-step FEMWARP was much
more robust.  This is probably because the reversals in the first
three rows are primarily Type 2 reversals since the meshes are very coarse,
whereas the ``tire'' mesh is finer and is therefore more likely to 
see Type 1 reversals according to the arguments given in the previous 
section.  Small-step FEMWARP is intended to fix Type 1 reversals but
is not effective against Type 2 reversals.

We also compared Opt-MS and 
the hybrid FEMWARP/Opt-MS method described in the
last section.  These results are shown in the fourth and fifth
columns of \tabref{tab:compare_3D_2}.  The minimum allowed stepsize
for the algorithms was again taken to be $\Delta \alpha$ as shown
in the last column of the table.
Opt-MS performed poorly because, as mentioned earlier,
it is not designed to handle large boundary motions.  The performance
of the hybrid is comparable, and in some cases superior, to that of 
small-step FEMWARP.


\begin{table}[!h]
\caption{Values of $\alpha^{\rm max}$ for FEMWARP on example three-dimensional
meshes.} 
\label{tab:compare_3D_1}
\begin{center}
\begin{tabular}{ccccc} \hline
Mesh name & \# bdry vertices & \# vertices & $\alpha^{\rm max}_{\rm FEMWARP}$ &
$\Delta\alpha$\\
\hline
foam5 & 1048 & 1337 & $0.7$ & $0.1$\\
gear & 606 & 866 &  $3.5$ & $0.1$ \\
hook & 790 & 1190 & $0.16$ & $0.01$ \\
tire & 1248 & 2570 & $0.15$ & $0.05$
\end{tabular}
\end{center}
\end{table}

\begin{table}[!h]
\caption{Values of $\alpha^{\rm max}$ for FEMWARP derivative algorithms on example three-dimensional
meshes.} 
\label{tab:compare_3D_2}
\begin{center}
\begin{tabular}{ccccccccc} \hline
Mesh name & $\alpha^{\rm max}_{\rm ssFEMWARP}$ & $\mbox{NCHOL}_{\rm SS}$&
$\alpha^{\rm max}_{\rm Opt-MS}$&
$\alpha^{\rm max}_{\rm hybrid}$&
$\Delta\alpha$\\
\hline
foam5 & $1.0$ & $4$ & $0$ & $1.2$ & $0.1$\\
gear & $3.5$ & $4$ & $0.6$ & $3.5$ & $0.1$ \\
hook & $0.16$ & $2$  & $0$ & $0.16$ & $0.01$ \\
tire & $1.60$ & $13$  & $0$ & $2.0$ & $0.05$
\end{tabular}
\end{center}
\end{table}

\section{Conclusions}
\label{sec:conclusions}

We studied an algorithm called FEMWARP for warping triangular and tetrahedral meshes.  
The first step in the algorithm is to determine a set of local weights for each 
interior vertex using finite element methods.  Second, a user-supplied deformation is 
applied to the boundary vertices.  The third and final step is to solve a system of linear 
equations based upon the weights and the new positions of the boundary vertices to determine 
the final positions of the interior vertices.

There are three main advantages of the FEMWARP algorithm as compared to
other mesh-updating methods.  First, if a continuous boundary
deformation is given, then FEMWARP is valid for computing the resulting
trajectory specifying the movement of the interior vertices.  In addition,
these trajectories will be continuous, which is vital for applications
where continuity of motion is required.  Second, sparse matrix algorithms
may be used to solve the linear system which determines the final
positions of the interior vertices.  Third, FEMWARP is exact if the boundary
deformation is affine.  

The main limitation of FEMWARP (as well as most other mesh-updating methods)
is that it can fail to yield valid deformed meshes, i.e., it sometimes
produces element reversals.  Element reversal was our main focus of the
paper.  We analyzed the case when FEMWARP produced element reversals 
and proposed four workarounds which include:  taking smaller steps, using a 
finer mesh, employing the mean value map to compute the weights, and, finally, using a hybrid algorithm 
which combines FEMWARP and Opt-MS. 

We tested the robustness of FEMWARP, small-step FEMWARP, a version of
FEMWARP which employed the mean value map to compute the weights, and hybrid
FEMWARP/Opt-MS on 2D annulus test cases and 3D general unstructured
meshes.  The use of the mean value map to compute the weights did not
improve the performance of FEMWARP on 2D meshes, and hence was not further 
considered.  Small-step FEMWARP and hybrid FEMWARP/Opt-MS generally
outperform plain FEMWARP, sometimes significantly.  However, an important limitation of the
hybrid FEM-WARP/Opt-MS algorithm is that there is no guarantee
that Opt-MS eliminates all element reversals.  In addition, the 
algorithm may not preserve the continuous trajectories needed by some
applications.     

Another limitation of the proposed techniques is that they do not
guarantee the quality of the warped mesh.  Future work should focus
on the development of mesh warping techniques with quality guarantees
for the resulting mesh.  For example, it may be possible to use 
optimization in conjunction with the proposed mesh warping
techniques in order to develop such quality guarantees.  Of course, 
the resulting algorithms will be likely be more expensive than the current 
techniques.

In order to reduce the time needed to compute the mesh warping, future
work will also focus on the development of an algorithm which formulates
a symmetric linear system (instead of forming the symmetric
positive definite linear system (\ref{eqn:lastsystem_framework})) expressing the coordinates
of each interior vertex in the original mesh in terms of its neighboring vertices.
Such a method would belong to the same more general class of methods to
which FEMWARP belongs.  This more general framework is called the
linear weighted Laplacian smoothing framework and was developed by
the first author in~\cite{ShontzThesis}.  Methods within this framework
share many of the same properties as FEMWARP.

The proposed mesh warping techniques in this paper are geometric in nature 
and do not take into account any knowledge of the particular PDE which may 
be causing the deformation.  Thus, another possibility for future work is 
to tie the mesh warping techniques to a particular PDE solver so that the 
deformed mesh is computed as a function of the PDE which is creating the motion in 
addition to being computed as a function of the particular domain geometry 
and the deformation upon it.  However, the use of such a technique is not 
always possible.  In particular, such a technique could not be used to
compute the deformed mesh in conjunction with a discrete set of motion data
stemming from a laboratory experiment.  The goal of this work was to 
develop improved mesh warping techniques which more robustly handle 
deformations in the cases where FEMWARP reverses elements.  Some of the
proposed techniques are applicable for problems with PDE-based motions, whereas
other techniques are applicable for problems with discrete datasets for 
mesh motions stemming from laboratory experiments.

\section{Acknowledgements}
\label{sec:acknowledgements}
The authors wish to
thank Dr.~Lori Diachin of Lawrence Livermore National Laboratory
and Dr.~Patrick Knupp of Sandia National Laboratories for providing us with 
the 3D test meshes.  They benefited from
conversations with G.~Bailey and A.~Schatz of Cornell, C.~Patron of Risk Capital, 
T.~Tomita of BRM Consulting, 
and especially H.~Kesten, formerly of Cornell.  
They also wish to thank the two anonymous referees for their careful reading 
of the paper and for their helpful suggestions which strengthened it.

\bibliographystyle{spmpsci}      
\bibliography{refs}   

\end{document}